\documentclass[11pt]{article}
\usepackage[centertags]{amsmath}
\usepackage{amsfonts}
\usepackage{amssymb}
\usepackage{amsthm}
\usepackage{newlfont}
\usepackage{graphicx}

\newfont{\bb}{msbm10 at 12pt}

\def\e{\hbox{\bf E}}
\def\t{\hbox{\bf T}}
\def\n{\hbox{\bf N}}
\def\b{\hbox{\bf B}}

\newtheorem{theorem}{Theorem}[section]
\newtheorem{definition}[theorem]{Definition}

\newtheorem{lemma}[theorem]{Lemma}
\newtheorem{example}[theorem]{Example}

\begin{document}
\title{Position vectors of a spacelike general helices in Minkowski Space $\e_1^3$}

\author{Ahmad T. Ali \\
Mathematics Department,\\
Faculty of Science, Al-Azhar University,\\
Nasr City, 11448, Cairo, Egypt.\\
E-mail: \textit{atali71@yahoo.com}}

\maketitle
\begin{abstract}
In this paper, position vector of a spacelike general helix with respect to standard frame in Minkowski space E$^3_1$ are studied in terms of Frenet equations. First, a vector differential equation of third order is constructed to determine position vector of an arbitrary spacelike general helix. In terms of solution, we determine the parametric representation of the general helices from the intrinsic equations. Moreover, we give some examples to illustrate how to find the position vectors of a spacelike general helices with a spacelike and timelike principal normal vector from the intrinsic equations.\\

\textbf{M.S.C. 2000}: 53C40, 53C50
\end{abstract}

\section{Introduction}
Natural scientists have long held a fascination, sometimes bordering on mystical obsession for helical structures in nature. Helices arise in nanosprings, carbon nanotubes, $\alpha-$helices, DNA double and collagen triple helix, the double helix shape is commonly associated with DNA, since the double helix is structure of DNA \cite{cam}. This fact was published for the first time by Watson and Crick in 1953 \cite{wats}. They constructed a molecular model of DNA in which there were two complementary, antiparallel (side-by-side in opposite directions) strands of the bases guanine, adenine, thymine and cytosine, covalently linked through phosphodiesterase bonds. Each strand forms a helix and two helices are held together through hydrogen bonds, ionic forces, hydrophobic interactions and van der Waals fores forming a double helix, lipid bilayers, bacterial flagella in Salmonella and E. coli, aerial hyphae in actynomycetes, bacterial shape in spirochetes, horns, tendrils, vines, screws, springs, helical staircases and sea shells (helico-spiral structures) \cite{chou, cook}.

In the local differential geometry, we think of curves as a geometric set of points, or locus. Intuitively, we are thinking of a curve as the path traced out by a particle moving in $\e^3$. So, investigating position vectors of the curves is a classical aim to determine behavior of the particle (or the curve, i.e.). There exists a vast literature on this subject, for instance \cite{ali2, chen, Il1, kula1, tur}.

A curve of constant slope or general helix is defined by the property that the tangent lines make a constant angle with a fixed direction. A necessary and sufficient condition that a curve to be general helix in Minkowski space $\e_1^3$ is that ratio of curvature to torsion be constant \cite{ferr}.

Indeed, a helix is a special case of the general helix, if both curvature $\kappa(s)$ and torsion $\tau(s)$ are non-zero constants, it is called a circular helix or only a $W$-curve \cite{Il2, Il3}. Its known that straight line ($\kappa(s)=0$) and circle ($\tau(s)=0$) are degenerate-helix examples \cite{kuhn, mill}. In fact, circular helix is the simplest three-dimensional spirals. One of the most interesting spiral example is $k$-Fibonacci spirals. These curves appear naturally from studying the $k$-Fibonacci numbers $\{F_{k,n}\}^{\infty}_{n=0}$ and the related hyperbolic $k$-Fibonacci function. Fibonacci numbers and the related Golden Mean or Golden section appear very often in theoretical physics and physics of the high energy particles \cite{elnasch1, elnasch2}. Three-dimensional $k$-Fibonacci spirals was studied from a geometric point of view in \cite{falc}.

All helices ($W-$curves) in $\e^3_1$ are completely classified by Walfare in \cite{walf}. For instance, the only planar space-like degenerate helices are circles and hyperbolas. In \cite{Il3}, the authors investigated position vectors of a time-like and a null helix ($W-$curve) with respect to Frenet frame.

For a given couple of one variable function (eventually curvature and torsion parameterized by arclength) one might like to get an arclength parameterized curve for which the couple works as the curvature and torsion functions. This problem is usually referred as "the solving natural equations problem". The natural equations for general helices can be integrated, nor only in Euclidean space $\e^3$, but also in 3-sphere $\mathbb{S}^3$ (the hyperbolic space is poor in this kind of curves and only helices are general helices). Indeed one uses the fact that general helices are geodesic either of right general cylinders or of Hopf cylinders, according to the curve lies in $\e^3$ or $\mathbb{S}^3$, respectively \cite{ali1, barr}.

In this work, we use vector differential equations established by means of Frenet equations in Minkowski space $\e^3_1$ to determine position vectors of the spacelike general helices according to standard frame in $\e^3_1$. We obtain the arclength parameterized of a space-like general helix with a spacelike and timelike principal normal vector from the intrinsic equations. We hope these results will be helpful to mathematicians who are specialized on mathematical modeling.

\section{Preliminaries}
To meet the requirements in the next sections, here, the basic elements
of the theory of curves in the space $\e_1^3$ are briefly presented (A
more complete elementary treatment can be found in \cite{Onei}.)\\

The Minkowski three-dimensional space $\e_1^3$ is the real vector space $\mathbb{R}^3$ endowed with the standard flat Lorentzian metric given by:
$$
g =-dx_1^2+dx_2^2+dx_3^2,
$$
where $(x_{1},x_{2},x_{3})$ is a rectangular coordinate system of $\e_1^3$. If $u=(u_1,u_2,u_3)$ and $v=(v_1,v_2,v_3)$ are arbitrary vectors in $\e_1^3$, we define the (Lorentzian) vector product of $u$ and $v$ as the following:
$$
u\times v=\Bigg|\begin{array}{ccc}
            i & -j & -k \\
            u_1 & u_2 & u_3 \\
            v_1 & v_2 & v_3
          \end{array}\Bigg|.
$$
Since $g$ is an indefinite metric, recall that a vector $v\in \e_1^3$ can have one of three Lorentzian characters: it can be space-like if $g(v,v)>0$ or $v=0$, timelike if $g(v,v)<0$ and null if $g(v,v)=0$ and $v\neq 0$. Similarly, an arbitrary curve $\psi =\psi(s)$ in $\e_1^3$ can locally be spacelike, timelike or null (lightlike), if all of its velocity vectors $\psi'$ are respectively spacelike, timelike or null (lightlike), for every $s\in I\subset \mathbb{R}$. The pseudo-norm of an arbitrary vector $a\in \e_1^3$ is given by $\left\Vert a\right\Vert=\sqrt{\left\vert g(a,a) \right\vert }$. $\psi$
is called an unit speed curve if velocity vector $v$ of $\psi$
satisfies $\left\Vert v\right\Vert =1$. For vectors $v,w\in\e_1^3$ it is said to be orthogonal if and only if $g(v,w) =0$.\\

Denote by $\left\{\t,\n,\b\right\}$ the moving Frenet frame along the curve $\psi$ in the space $\e_1^3$. For an arbitrary curve $\psi$ with first and second curvature, $\kappa$ and $\tau$ in the space $\e_1^3$, the following Frenet formulae are given in \cite{Il3}:\\

If $\psi$ is a spacelike curve with a spacelike or timelike principal normal vector $\n$, then the Frenet formulae read
\begin{equation}\label{u1}
\left[
\begin{array}{c}
\t'\\
\n'\\
\b'\\
\end{array}
\right]
=\left[
  \begin{array}{ccc}
    0 & \kappa  & 0 \\
-\varepsilon\kappa  & 0 & \tau  \\
0 & \tau  & 0\\
  \end{array}
\right]
\left[
\begin{array}{c}
\t \\
\n \\
\b\\
\end{array}
\right],
\end{equation}
where
$$
g(\t,\t)=1,\,\,g(\n,\n)=\varepsilon=\pm1,\,\,g(\b,\b)=-\varepsilon,\,\,g(\t,\n)=g(\t,\b)=g(\n,\b)=0.
$$
The angle between two vectors in Minkowski space is defined by \cite{bili}:
\begin{definition}\label{def-1}
Let $X$ and $Y$ be spacelike vectors in $\e_1^3$ that span a spacelike vector subspace, then we have $|g(X,Y)|\leq\|X\|\|Y\|$ and hence, there is a unique positive real number $\phi$ such that
$$
|g(X,Y)|=\|X\|\|Y\|\mathrm{cos}[\phi].
$$
The real number $\phi$ is called the Lorentzian spacelike angle between $X$ and $Y$.
\end{definition}

\begin{definition}\label{def-2}
Let $X$ and $Y$ be spacelike vectors in $\e_1^3$ that span a timelike vector subspace, then we have $|g(X,Y)|>\|X\|\|Y\|$ and hence, there is a unique positive real number $\phi$ such that
$$
|g(X,Y)|=\|X\|\|Y\|\mathrm{cosh}[\phi].
$$
The real number $\phi$ is called the Lorentzian timelike angle between $X$ and $Y$.
\end{definition}

\begin{definition}\label{def-3}
Let $X$ be a spacelike vector and $Y$ a positive timelike vector in $\e_1^3$, then there is a unique non-negative real number $\phi$ such that
$$
|g(X,Y)|=\|X\|\|Y\|\mathrm{sinh}[\phi].
$$
The real number $\phi$ is called the Lorentzian timelike angle between $X$ and $Y$.
\end{definition}

\begin{definition}\label{def-4}
Let $X$ and $Y$ be positive (negative) timelike vectors in $\e_1^3$, then there is a unique non-negative real number $\phi$ such that
$$
g(X,Y)=\|X\|\|Y\|\mathrm{cosh}[\phi].
$$
The real number $\phi$ is called the Lorentzian timelike angle between $X$ and $Y$.
\end{definition}

\section{Position vector of a spacelike general helix}
The problem of the determination of parametric representation of the position vector of an arbitrary space curve according to the intrinsic equations is still open in the Euclidean space E$^3$ and in the Minkowski spae E$_1^3$ \cite{ali1, eisenh, lips}. This problem is not easy to solve in general case. However, this problem is solved in three special cases only, Firstly, in the case of a plane curve ($\tau=0$). Secondly, in the case of a helix ($\kappa$ and $\tau$ are both non-vanishing constant). Recently, Ali and Turgut \cite{ali3} adapted fundamental existence and uniqueness theorem for space curves to time-like curves of the space $\e_1^3$ and constructed a vector differential equation to solve this problem in the case of a general helix ($\frac{\tau}{\kappa}$ is constant) in Minkowski space $\e_1^3$. However, this problem is not solved in other cases of the space curve.\\

Our main result in this work is to determine the parametric representation of the position vector $\psi$ from intrinsic equations in E$_1^3$ for the case of a spacelike general helix with a spacelike and timelike principal normal vector. In the light of our main problem, first we give:

\begin{theorem}\label{th-main} Let $\psi=\psi(s)$ be an unit speed spacelike curve parameterized by the arclength $s$. Suppose $\psi=\psi(\theta)$ is another parametric representation of this curve by the parameter $\theta=\int\kappa(s)ds$.  Then, the tangent vector $\t$ satisfies a vector differential equation of third order as follows:
\begin{equation}\label{u2}
\Big(\frac{1}{f(\theta)}\t''(\theta)\Big)'+\Big(\frac{\varepsilon-f^2(\theta)}{f(\theta)}\Big)\t'(\theta)
-\varepsilon\frac{f'(\theta)}{f(\theta)}\t(\theta)=0,
\end{equation}
where $f(\theta)=\frac{\tau(\theta)}{\kappa(\theta)}$.
\end{theorem}

{\bf Proof.} Let $\psi=\psi(s)$ be an unit speed spacelike curve. If we write this curve in the another parametric representation $\psi=\psi(\theta)$, where $\theta=\int\kappa(s)ds$, we have new Frenet equations as follows:
\begin{equation}\label{u3}
 \left[
   \begin{array}{c}
     \t'(\theta) \\
     \n'(\theta) \\
     \b'(\theta) \\
   \end{array}
 \right]=\left[
           \begin{array}{ccc}
             0 & 1 & 0 \\
             -\varepsilon & 0 & f(\theta) \\
             0 & f(\theta) & 0 \\
           \end{array}
         \right]\left[
   \begin{array}{c}
     \t(\theta) \\
     \n(\theta) \\
     \b(\theta) \\
   \end{array}
 \right],
 \end{equation}
where $f(\theta)=\frac{\tau(\theta)}{\kappa(\theta)}$. If we substitute the first equation of the new Frenet equations (\ref{u3}) to second equation of (\ref{u3}), we have
\begin{equation}\label{u4}
\b(\theta)=\frac{1}{f(\theta)}\Big[\t''(\theta)+\varepsilon\t(\theta)\Big].
\end{equation}
Substituting the above equation in the last equation from (\ref{u3}), we obtain a vector differential equation of third order (\ref{u2}) as desired.

The equation (\ref{u2}) is not easy to solve in general case. If one solves this equation, the natural representation of the position vector of an arbitrary space curve can be determined in natural representation form as follows:
\begin{equation}\label{u5}
\psi(s)=\int\,\t(s)\,ds+C,
\end{equation}
or in parametric representation form
\begin{equation}\label{u6}
\psi(\theta)=\int\frac{1}{\kappa(\theta)}\t(\theta)d\theta+C,
\end{equation}
where $\theta=\int\kappa(s)ds$.

We can solve the equation (\ref{u2}) in the case of a spacelike general helix. Because the tangent vector $\t$ is a spacelike vector for the spacelike curve, we can give the following three lemmas using the definitions \ref{def-1}, \ref{def-2} and \ref{def-3}. The following propositions are \textit{new characterizations} for spacelike general helices in $\e_1^3$:

\begin{lemma}\label{lm-1}
Let $\psi:I\rightarrow\e_1^3$ be a spacelike curve that is parameterized by arclength with intrinsic equations $\kappa=\kappa(s)$ and $\tau=\tau(s)$. The curve $\psi$ is a general helix (its tangent vectors make a constant Lorentzian spacelike angle $\phi=\pm\mathrm{arccos}[n]$, with a fixed spacelike straight line in the space) if and only if $\psi$ is a spacelike curve with a timelike principal normal vector and $\frac{\tau(s)}{\kappa(s)}=\mp\cot[\phi]$.
\end{lemma}

{\bf Proof:} $(\Leftarrow)$ Let $\mathbf{d}$ be the unitary fixed spacelike vector makes a constant spacelike angle, $\phi=\pm\mathrm{arccos}[n]$, with the tangent vector $\t$. Therefore
\begin{equation}\label{u7}
g(\mathbf{d},\t)=n.
\end{equation}
Differentiating the equation (\ref{u7}) with respect to the variable $s$ and using the Frenet frame (\ref{u1}), we get
\begin{equation}\label{u8}
\kappa\,g(\mathbf{d},\n)=0.
\end{equation}
The curvature $\kappa(s)$ do not equal zero, therefore the vector $\mathbf{d}$ is orthogonal to $\n$ and take the form
\begin{equation}\label{u9}
\mathbf{d}=\cos[\phi]\t+\lambda\,\b.
\end{equation}
From the unitary of the vector $\mathbf{d}$ we get
\begin{equation}\label{u10}
\lambda^2=\frac{\sin^2[\phi]}{-\varepsilon}.
\end{equation}
Therefore $\varepsilon$ must equal $-1$, $\lambda=\pm\sin[\phi]$, and hence $\psi$ is a spacelike curve with a timelike principal normal vector. If we differentiate the equation (\ref{u9}), we obtain the desired result.

$(\Leftarrow)$ Suppose that $\psi$ is a spacelike curve with a timelike principal normal vector and $\tau(s)=\mp\cot[\phi]\kappa(s)$. Let us consider the vector
$$
\mathbf{d}=\cos[\phi]\t\pm\sin[\phi]\b.
$$
We will prove that the vector $\mathbf{d}$ is a constant vector. Indeed, applying Frenet formula (\ref{u1})
$$
\mathbf{d}'(s)=\Big(\cos[\phi]\kappa(s)\pm\sin[\phi]\tau(s)\Big)\n=0.
$$
Therefore, the vector $\mathbf{d}$ is constant. This concludes the proof of lemma (\ref{lm-1}).

\begin{lemma}\label{lm-2}
Let $\psi:I\rightarrow\e_1^3$ be a spacelike curve that is parameterized by arclength with intrinsic equations $\kappa=\kappa(s)$ and $\tau=\tau(s)$. The curve $\psi$ is a general helix (its tangent vectors make a constant Lorentzian timelike angle $\phi=\pm\mathrm{arccosh}[n]$, with a fixed spacelike straight line in the space) if and only if $\psi$ is a spacelike curve with a spacelike principal normal vector and $\Big|\frac{\tau(s)}{\kappa(s)}\Big|=\Big|\coth[\phi]\Big|>1$.
\end{lemma}

\begin{lemma}\label{lm-3}
Let $\psi:I\rightarrow\e_1^3$ be a spacelike curve that is parameterized by arclength with intrinsic equations $\kappa=\kappa(s)$ and $\tau=\tau(s)$. The curve $\psi$ is a general helix (its tangent vectors make a constant Lorentzian timelike angle $\phi=\pm\mathrm{arccosh}[n]$, with a fixed timelike straight line in the space) if and only if $\psi$ is a spacelike curve with a spacelike principal normal vector and $\Big|\frac{\tau(s)}{\kappa(s)}\Big|=\Big|\tanh[\phi]\Big|<1$.
\end{lemma}

The proof of the lemmas (\ref{lm-2}) and $(\ref{lm-3})$ are similar as the proof of the lemma (\ref{lm-1}). Now we give three theorems correspondence to the three lemmas as follows:

\begin{theorem}\label{th-main2} The position vector $\psi$ of a spacelike general helix with a timelike principal normal vector, whose tangent vector makes a constant Lorentzian spacelike angle, with a fixed spacelike straight line in the space, is computed in the natural representation form:
\begin{equation}\label{u14}
\psi(s)=\sqrt{1-n^2}\int\Big(\sinh\Big[\sqrt{1+m^2}\,\int\kappa(s)ds\Big],
\cosh\Big[\sqrt{1+m^2}\,\int\kappa(s)ds\Big],m\Big)ds,
\end{equation}
or in the parametric form
\begin{equation}\label{u15}
\psi(\theta)=\int\frac{\sqrt{1-n^2}}{\kappa(\theta)}\Big(\sinh[\sqrt{1+m^2}\,\theta],
\cosh[\sqrt{1+m^2}\,\theta],m\Big)d\theta,
\end{equation}
where $\theta=\int\kappa(s)ds$, $m=\frac{n}{\sqrt{1-n^2}}$, $n=\cos[\phi]$ and $\phi$ is the spacelike angle between the fixed spacelike straight line $\mathbf{e}_3$ (axis of a spacelike general helix with a timelike principal normal vector) and the tangent vector of the curve.
\end{theorem}

{\bf Proof:} If $\psi(\theta)$ is a spacelike general helix with the timelike principal normal vector whose tangent vector $\t$ makes a spacelike angle $\phi=\pm\mathrm{\arccos}[n]$ with a straight spacelike line $U$, then we can write $\varepsilon=-1$, $f(\theta)=\cot[\phi]=m$, where $\theta=\int\kappa(s)ds$ and $m=\frac{n}{\sqrt{1-n^2}}$. Therefore the equation (\ref{u2}) becomes
\begin{equation}\label{u17}
\t'''(\theta)-(1+m^2)\t'(\theta)=0.
\end{equation}
If we write the tangent vector as the following:
\begin{equation}\label{u18}
\t=T_1(\theta)\mathbf{e}_1+T_2(\theta)\mathbf{e}_2+T_3(\theta)\mathbf{e}_3.
\end{equation}
Now, the curve $\psi$ is a spacelike general helix, i.e. the tangent vector $\t$ makes a constant spacelike angle, $\phi$, with the constant spacelike vector called the axis of the general helix. So, with out loss of generality, we take the axis of a general helix is parallel to the spacelike vector $\bold{e}_3$. Then
\begin{equation}\label{u19}
T_3=g(\t,\mathbf{e}_3)=n.
\end{equation}
On other hand the tangent vector $\t$ is a unit spacelike vector, so the following condition is satisfied
\begin{equation}\label{u20}
-T_1^2(\theta)+T_2^2(\theta)=1-n^2.
\end{equation}
The general solution of equation (\ref{u20}) can be written in the following form:
\begin{equation}\label{u21}
\begin{array}{ll}
T_1=\sqrt{1-n^2}\sinh[t(\theta)],\,\,\,\,\,T_2=\sqrt{1-n^2}\cosh[t(\theta)],
\end{array}
\end{equation}
where $t$ is an arbitrary function of $\theta$. Every component of the vector $\t$ is satisfied the equation (\ref{u17}). So, substituting the components $T_1(\theta)$ and $T_2(\theta)$ in the equation (\ref{u17}), we have the following differential equations of the function $t(\theta)$
\begin{equation}\label{u22}
3t't''\sinh[t]-\Big[(1+m^2)t'-t'^3-t'''\Big]\cosh[t]=0,
\end{equation}
\begin{equation}\label{u23}
3t't''\cosh[t]-\Big[(1+m^2)t'-t'^3-t'''\Big]\sinh[t]=0.
\end{equation}
It is easy to prove that the above two equations lead to the following two equations:
\begin{equation}\label{u24}
t'\,t''=0,
\end{equation}
\begin{equation}\label{u25}
(1+m^2)t'-t'^3-t'''=0.
\end{equation}
Because, the parameter $t$ is a variable (not constant), then $t'\neq0$, so that the general solution of the equation (\ref{u24}) is
\begin{equation}\label{u26}
t(\theta)=c_2+c_1\,\theta,
\end{equation}
where $c_1$ and $c_2$ are constants of integration. The constant $c_2$ can be disappear if we change the parameter $t\rightarrow t+c_2$. Substituting the solution (\ref{u26}) in the equation (\ref{u25}), we obtain the following condition:
$$
c_1^2-1-m^2=0
$$
which leads to $c_1=\sqrt{1+m^2}$.

Now, the tangent vector take the following form:
\begin{equation}\label{u27}
\t(\theta)=\sqrt{1-n^2}\Big(\sinh[\sqrt{1+m^2}\,\theta],
\cosh[\sqrt{1+m^2}\,\theta],m\Big).
\end{equation}
If we substitute the equation (\ref{u27}) in the two equations (\ref{u5}) and (\ref{u6}), we have the two equations (\ref{u14}) and (\ref{u15}), which it completes the proof.

\begin{theorem}\label{th-main3} The position vector $\psi$ of a spacelike general helix with a spacelike principal normal vector, whose tangent vector makes a constant Lorentzian timelike angle, with a fixed spacelike straight line in the space, is computed in the natural representation form:
\begin{equation}\label{u141}
\psi(s)=\sqrt{n^2-1}\int\Big(\cosh\Big[\sqrt{m^2-1}\,\int\kappa(s)ds\Big],
\sinh\Big[\sqrt{m^2-1}\,\int\kappa(s)ds\Big],m\Big)ds,
\end{equation}
or in the parametric form
\begin{equation}\label{u151}
\psi(\theta)=\int\frac{\sqrt{n^2-1}}{\kappa(\theta)}\Big(\cosh[\sqrt{m^2-1}\,\theta],
\sinh[\sqrt{m^2-1}\,\theta],m\Big)d\theta,
\end{equation}
where $\theta=\int\kappa(s)ds$, $m=\frac{n}{\sqrt{n^2-1}}$, $n=\cosh[\phi]$ and $\phi$ is the timelike angle between the fixed spacelike straight line $\mathbf{e}_3$ (axis of a spacelike general helix with a spacelike principal normal vector) and the tangent vector of the curve.
\end{theorem}

\begin{theorem}\label{th-main4} The position vector $\psi$ of a spacelike general helix with a spacelike principal normal vector, whose tangent vector makes a constant Lorentzian timelike angle, with a fixed timelike straight line in the space, is computed in the natural representation form:
\begin{equation}\label{u142}
\psi(s)=\sqrt{n^2+1}\int\Big(m,\cos\Big[\sqrt{1-m^2}\,\int\kappa(s)ds\Big],
\sin\Big[\sqrt{1-m^2}\,\int\kappa(s)ds\Big]\Big)ds,
\end{equation}
or in the parametric form
\begin{equation}\label{u152}
\psi(\theta)=\int\frac{\sqrt{n^2+1}}{\kappa(\theta)}\Big(m,\cos[\sqrt{1-m^2}\,\theta],
\sin[\sqrt{1-m^2}\,\theta]\Big)d\theta,
\end{equation}
where $\theta=\int\kappa(s)ds$, $m=\frac{n}{\sqrt{1+n^2}}$, $n=\sinh[\phi]$ and $\phi$ is the timelike angle between the fixed timelike straight line $\mathbf{e}_3$ (axis of a spacelike general helix with a spacelike principal normal vector) and the tangent vector of the curve.
\end{theorem}

According to lemmas \ref{lm-2} and \ref{lm-3}, the prove of the two theorems \ref{th-main3} and \ref{th-main4} are similar as the proof of the theorem \ref{th-main2}.

\section{Examples}

In this section, we take several choices for the curvature $\kappa$ and torsion $\tau$, and next, we apply Theorem \ref{th-main2}.

\begin{example}
The case of a spacelike general helix with
\begin{equation}\label{u28}
\kappa=\kappa(s),\,\,\,\,\,\tau=0,
\end{equation}
which is the intrinsic equations of a spacelike plane curve. There is a three cases correspondence to three theorems \ref{th-main2}, \ref{th-main3} and \ref{th-main4}.
\end{example}

{\bf Case 1:} If $\psi$ is a spacelike plane curve with a timelike principal normal vector then the tangent vectors make a constant Lorentzian spacelike angle $\phi=\mathrm{arccos}[n]$ with a fixed spacelike straight line $\mathbf{e}_3$. According to lemma \ref{lm-1}, we have $\cot[\phi]=\frac{\tau(s)}{\kappa(s)}=0$ which leads to $n=\cos[\phi]=0$ and $m=\cot[\phi]=0$. Substituting values $(n=m=0)$ in the equations (\ref{u14}) and (\ref{u15}) we have the explicit natural and parametric representation of such curve as follows:
\begin{equation}\label{u143}
\psi(s)=\int\Big(\sinh\Big[\int\kappa(s)ds\Big],
\cosh\Big[\int\kappa(s)ds\Big],0\Big)ds,
\end{equation}
\begin{equation}\label{u153}
\psi(\theta)=\int\frac{1}{\kappa(\theta)}\Big(\sinh[\theta],
\cosh[\theta],0\Big)d\theta,
\end{equation}
where $\theta=\int\kappa(s)ds$. Now, we give the parametric representation of a special example of a plane curve.

The position vector $\psi$ of a spacelike plane curve ($\kappa(s)=\frac{a}{a^2-s^2}$) with a timelike principal normal vector, whose tangent vector makes a constant Lorentzian spacelike angle, $\phi$, with a fixed spacelike straight line in the space, takes the form:
\begin{equation}\label{u1451}
\psi(s)=a\Big(-\mathrm{sech}[\theta],2\,\mathrm{arctan}\Big[\tanh[\frac{\theta}{2}]\Big],0\Big),
\end{equation}
where $s=a\,\tanh[\theta]$ and $\kappa(\theta)=\frac{\cosh^2[\theta]}{a}$. One can see a special example of such curve when $a=2$ in the lift hand side of the figure 1.

{\bf Case 2:} If $\psi$ is a spacelike plane curve with a spacelike principal normal and the tangent vectors make a constant Lorentzian timelike angle $\phi=\mathrm{arccosh}[n]$ with a fixed spacelike straight line $\mathbf{e}_3$. According to lemma \ref{lm-2}, we have $\coth[\phi]=\frac{\tau(s)}{\kappa(s)}=0$ which is contradiction $(\coth[\phi]\neq0)$. Therefore, we can write the following lemma.

\begin{lemma}\label{lm-4}
There is no a spacelike plane curve with a spacelike principal normal and the tangent vector makes a constant Lorentzian timelike angle with a fixed spacelike straight line.
\end{lemma}

{\bf Case 3:} If $\psi$ is a spacelike plane curve with a spacelike principal normal and the tangent vectors make a constant Lorentzian timelike angle $\phi=\mathrm{arcsinh}[n]$ with a fixed timelike straight line $\mathbf{e}_1$. According to lemma \ref{lm-3}, we have $\tanh[\phi]=\frac{\tau(s)}{\kappa(s)}=0$ which leads to $n=\sinh[\phi]=0$ and $m=\tanh[\phi]=0$. Substituting values $(n=m=0)$ in the equations (\ref{u142}) and (\ref{u152}) we have the explicit natural and parametric representation of such curve as follows:
\begin{equation}\label{u144}
\psi(s)=\int\Big(0,\cos\Big[\int\kappa(s)ds\Big],
\sin\Big[\int\kappa(s)ds\Big]\Big)ds,
\end{equation}
or in the parametric form
\begin{equation}\label{u154}
\psi(\theta)=\int\frac{1}{\kappa(\theta)}\Big(0,\cos[\theta],
\sin[\theta]\Big)d\theta,
\end{equation}
where $\theta=\int\kappa(s)ds$. Now, we give the parametric representation of a special example of a plane curve.

The position vector $\psi$ of a spacelike plane curve ($\kappa(s)=\frac{a}{a^2+s^2}$) with a space principal normal vector, whose tangent vector makes a constant Lorentzian timelike angle, $\phi$, with a fixed timelike straight line in the space, takes the form:
\begin{equation}\label{u1451}
\psi(s)=a\Big(0,2\,\mathrm{arctanh}\Big[\tan[\frac{\theta}{2}]\Big],\mathrm{sec}[\theta]\Big),
\end{equation}
where $s=a\,\tan[\theta]$ and $\kappa(\theta)=\frac{\cos^2[\theta]}{a}$. One can see a special example of such curve when $a=\frac{1}{2}$ in the right hand side of the figure 1.

\begin{figure}[ht]
\begin{center}
\includegraphics[width=4.5cm]{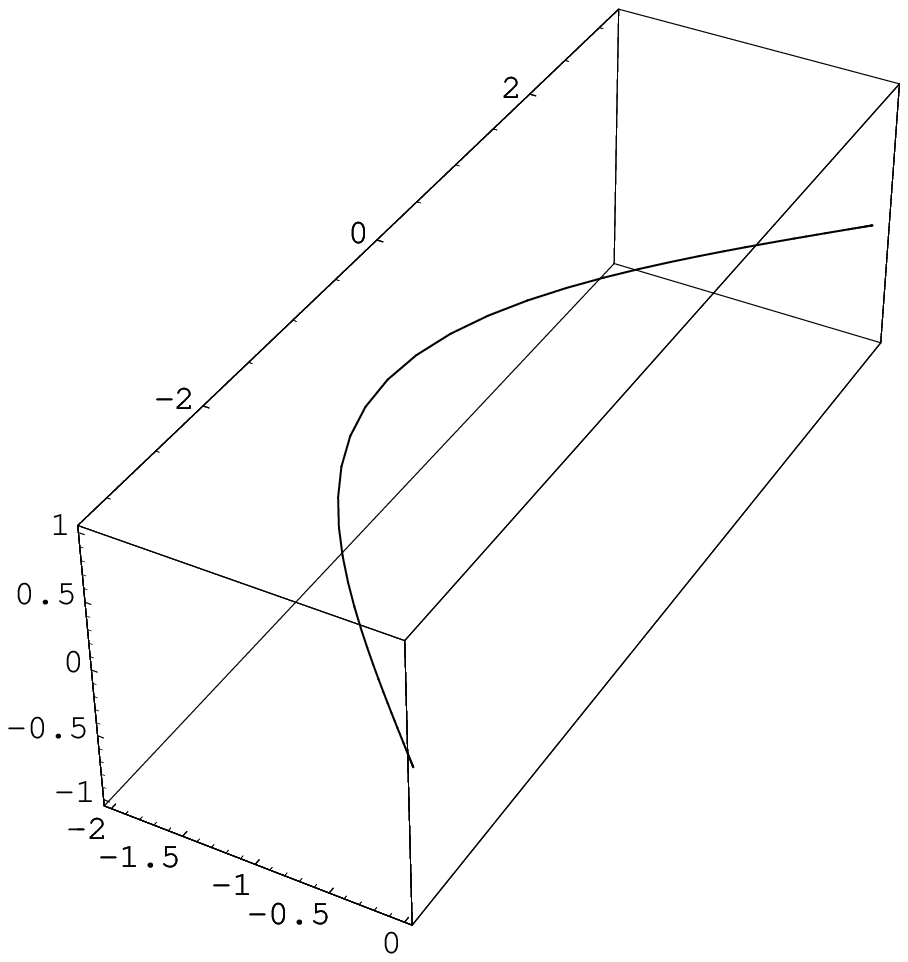}
\hspace*{0.5cm}
\includegraphics[width=4.0cm]{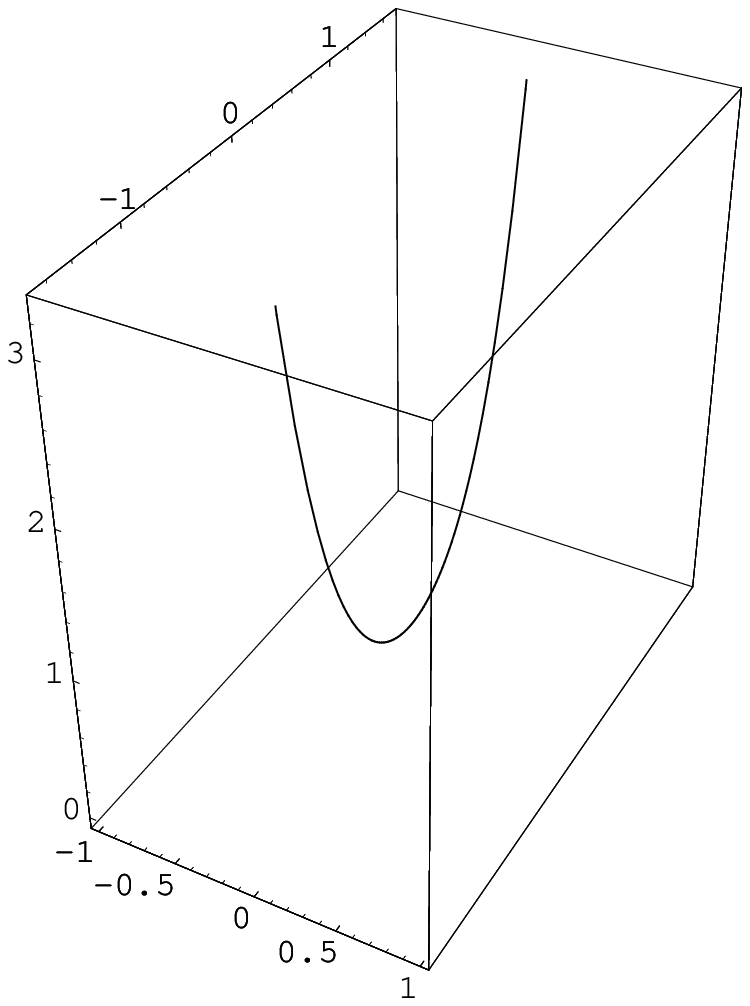}
\end{center}\caption{Some spacelike plane curves.}
\label{fig-x1}
\end{figure}

\begin{example}
The case of a spacelike general helix with
\begin{equation}\label{u28}
\kappa=\mathrm{constant},\,\,\,\,\,\tau=\mathrm{constant},
\end{equation}
which is the intrinsic equations of a spacelike W-curve or helix. There is a three cases correspondence to three theorems \ref{th-main2}, \ref{th-main3} and \ref{th-main4}.
\end{example}

{\bf Case 1:} The position vector $\psi$ of a spacelike W-curve with a timelike principal normal vector, whose tangent vector makes a constant Lorentzian spacelike angle, $\phi$, with a fixed spacelike straight line in the space, is computed in the parametric representation form:
\begin{equation}\label{u145}
\psi(s)=\frac{\kappa}{\kappa^2+\tau^2}\Big(\cosh[\xi],
\sinh[\xi],\frac{\tau}{\kappa}\,\xi\Big),
\end{equation}
where $\xi=\sqrt{\kappa^2+\tau^2}\,s$ and $\phi=\mathrm{arccot}[\frac{\tau}{\kappa}]$. One can see a special example of such curve when $\kappa=3$ and $\tau=2$ in the lift hand side of the figure 2.

{\bf Case 2:} The position vector $\psi$ of a spacelike W-curve with a spacelike principal normal vector, whose tangent vector makes a constant Lorentzian timelike angle, $\phi$, with a fixed spacelike straight line in the space, is computed in the parametric representation form:
\begin{equation}\label{u146}
\psi(s)=\frac{\kappa}{\tau^2-\kappa^2}\Big(\sinh[\xi],
\cosh[\xi],\frac{\tau}{\kappa}\,\xi\Big),
\end{equation}
where $\xi=\sqrt{\tau^2-\kappa^2}\,s$ and $\phi=\mathrm{arccoth}[\frac{\tau}{\kappa}]$. One can see a special example of such curve when $\kappa=1$ and $\tau=2$ in the middle of the figure 2.

{\bf Case 3:} The position vector $\psi$ of a spacelike W-curve with a spacelike principal normal vector, whose tangent vector makes a constant Lorentzian timelike angle, $\phi$, with a fixed timelike straight line in the space, is computed in the parametric representation form:
\begin{equation}\label{u1461}
\psi(s)=\frac{\kappa}{\kappa^2-\tau^2}\Big(\frac{\tau}{\kappa}\,\xi,\sin[\xi],
\cos[\xi]\Big),
\end{equation}
where $\xi=\sqrt{\kappa^2-\tau^2}\,s$ and $\phi=\mathrm{arctanh}[\frac{\tau}{\kappa}]$. One can see a special example of such curve when $\kappa=2$ and $\tau=1$ in the right hand side of the figure 2.

\begin{figure}[ht]
\begin{center}
\includegraphics[width=4.0cm]{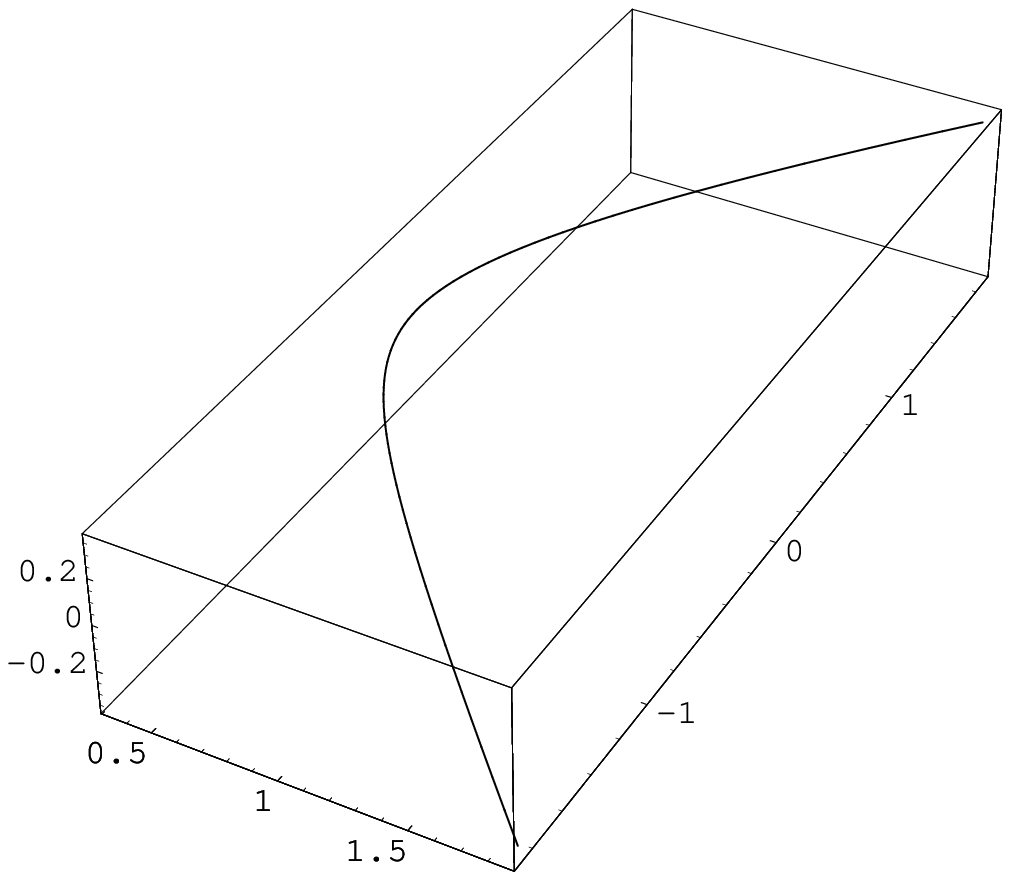}
\hspace*{0.5cm}
\includegraphics[width=3.5cm]{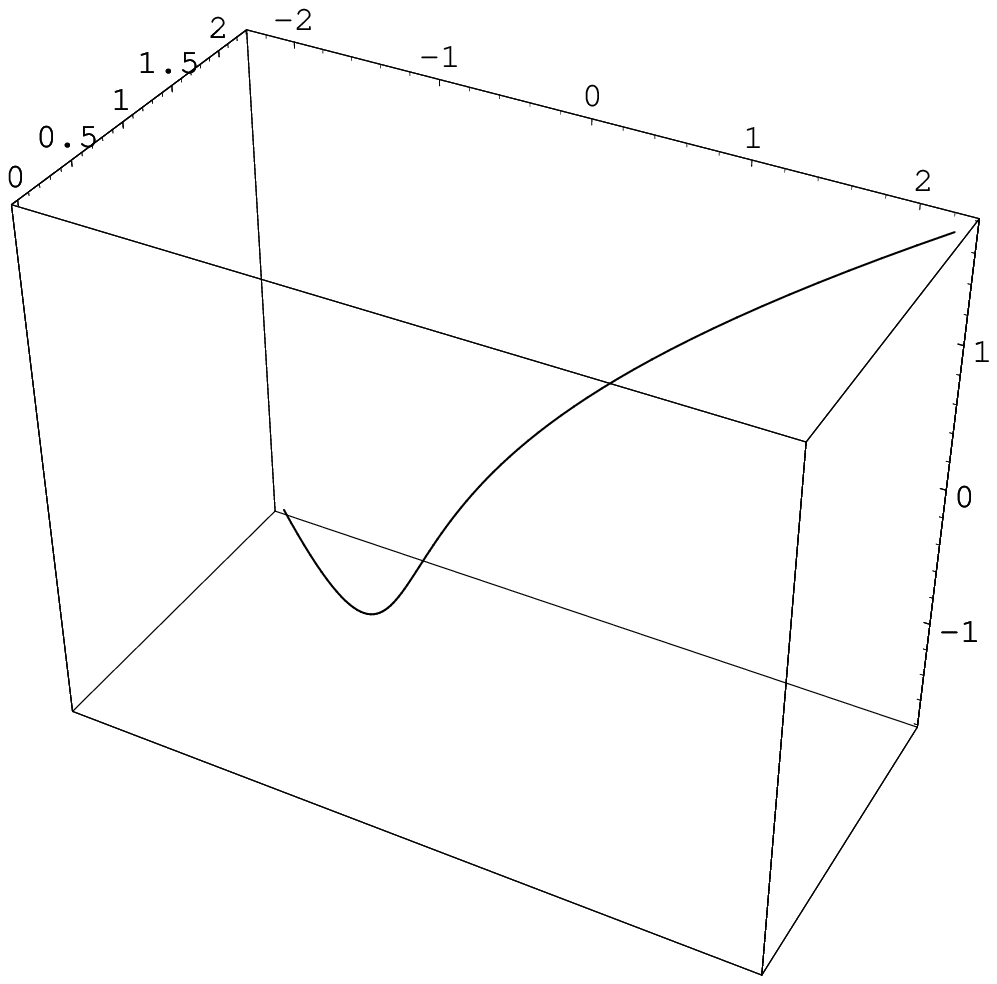}
\hspace*{0.5cm}
\includegraphics[width=5.5cm]{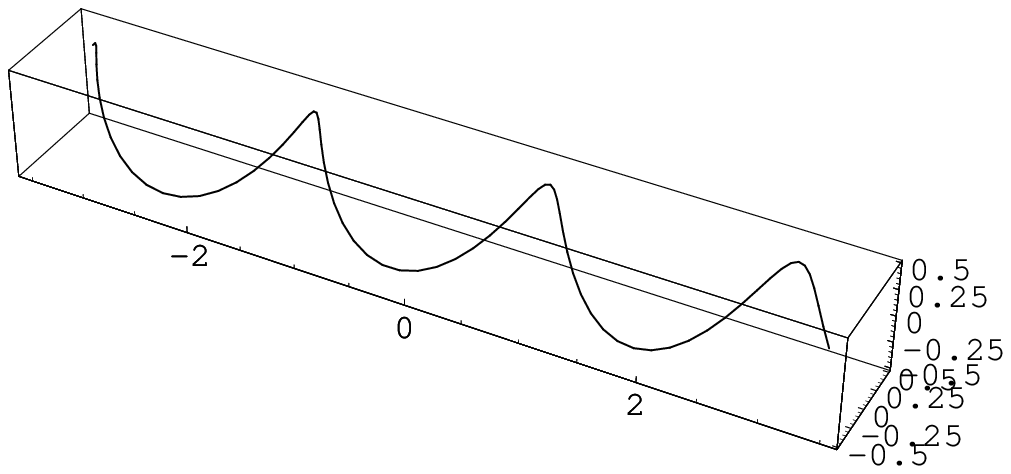}
\end{center}\caption{Some spacelike W-curves.}
\label{fig-x2}
\end{figure}

\begin{example}
The case of a spacelike general helix with
\begin{equation}\label{u28}
\kappa=\frac{h}{s},\,\,\,\,\,\tau=\frac{r}{s},
\end{equation}
where $h$ and $r$ are arbitrary constants. There is a three cases correspondence to three theorems \ref{th-main2}, \ref{th-main3} and \ref{th-main4}.
\end{example}

{\bf Case 1:} The position vector $\psi=(\psi_1,\psi_2,\psi_3)$ of a spacelike general helix with a timelike principal normal vector, whose tangent vector makes a constant Lorentzian spacelike angle, $\phi$, with a fixed spacelike straight line in the space, is computed in the parametric representation form:
\begin{equation}\label{u148}
\left\{
  \begin{array}{ll}
    \psi_1(\theta)=\frac{h\,\mathrm{e}^{\theta/h}}{h^2+r^2-1}\Big(\cosh\Big[\frac{\sqrt{h^2+r^2}}{h}\,\theta\Big]-
\frac{1}{\sqrt{h^2+r^2}}\sinh\Big[\frac{\sqrt{h^2+r^2}}{h}\,\theta\Big]\Big),\\
    \psi_2(\theta)=\frac{h\,\mathrm{e}^{\theta/h}}{h^2+r^2-1}\Big(\sinh\Big[\frac{\sqrt{h^2+r^2}}{h}\,\theta\Big]-
\frac{1}{\sqrt{h^2+r^2}}\cosh\Big[\frac{\sqrt{h^2+r^2}}{h}\,\theta\Big]\Big),\\
    \psi_3(\theta)=\frac{r\,\mathrm{e}^{\theta/h}}{\sqrt{h^2+r^2}},
  \end{array}
\right.
\end{equation}
where $s=\mathrm{e}^{\theta/h}$ and $\phi=\mathrm{arccot}[\frac{r}{h}]$. One can see a special example of such curve when $h=2$ and $r=1$ in the lift hand side of figure 3.

{\bf Case 2:} The position vector $\psi$ of a spacelike general helix with a spacelike principal normal vector, whose tangent vector makes a constant Lorentzian timelike angle, $\phi$, with a fixed spacelike straight line in the space, is computed in the parametric representation form:
\begin{equation}\label{u1481}
\left\{
  \begin{array}{ll}
    \psi_1(\theta)=\frac{h\,\mathrm{e}^{\theta/h}}{1+h^2-r^2}\Big(
\frac{1}{\sqrt{r^2-h^2}}\cosh\Big[\frac{\sqrt{r^2-h^2}}{h}\,\theta\Big]-\sinh\Big[\frac{\sqrt{r^2-h^2}}{h}\,\theta\Big]\Big),\\
    \psi_2(\theta)=\frac{h\,\mathrm{e}^{\theta/h}}{1+h^2-r^2}\Big(
\frac{1}{\sqrt{r^2-h^2}}\sinh\Big[\frac{\sqrt{r^2-h^2}}{h}\,\theta\Big]-\cosh\Big[\frac{\sqrt{r^2-h^2}}{h}\,\theta\Big]\Big),\\
    \psi_3(\theta)=\frac{r\,\mathrm{e}^{\theta/h}}{\sqrt{r^2-h^2}},
  \end{array}
\right.
\end{equation}
where $s=\mathrm{e}^{\theta/h}$ and $\phi=\mathrm{arccoth}[\frac{r}{h}]$. One can see a special example of such curve when $h=1$ and $r=4$ in the middle of the figure 2.

{\bf Case 3:} The position vector $\psi$ of a spacelike general helix with a spacelike principal normal vector, whose tangent vector makes a constant Lorentzian timelike angle, $\phi$, with a fixed timelike straight line in the space, is computed in the parametric representation form:
\begin{equation}\label{u14811}
\left\{
  \begin{array}{ll}
    \psi_1(\theta)=\frac{r\,\mathrm{e}^{\theta/h}}{\sqrt{h^2-r^2}},\\
\psi_2(\theta)=\frac{h\,\mathrm{e}^{\theta/h}}{1+h^2-r^2}\Big(
\frac{1}{\sqrt{h^2-r^2}}\cos\Big[\frac{\sqrt{h^2-r^2}}{h}\,\theta\Big]+\sin\Big[\frac{\sqrt{h^2-r^2}}{h}\,\theta\Big]\Big),\\
    \psi_2(\theta)=\frac{h\,\mathrm{e}^{\theta/h}}{1+h^2-r^2}\Big(
\frac{1}{\sqrt{h^2-r^2}}\sin\Big[\frac{\sqrt{h^2-r^2}}{h}\,\theta\Big]-\cos\Big[\frac{\sqrt{h^2-r^2}}{h}\,\theta\Big]\Big),\\
      \end{array}
\right.
\end{equation}
where $s=\mathrm{e}^{\theta/h}$ and $\phi=\mathrm{arctanh}[\frac{r}{h}]$. One can see a special example of such curve when $h=6$ and $r=1$ in the right hand side of figure 2.

\begin{figure}[ht]
\begin{center}
\includegraphics[width=5.0cm]{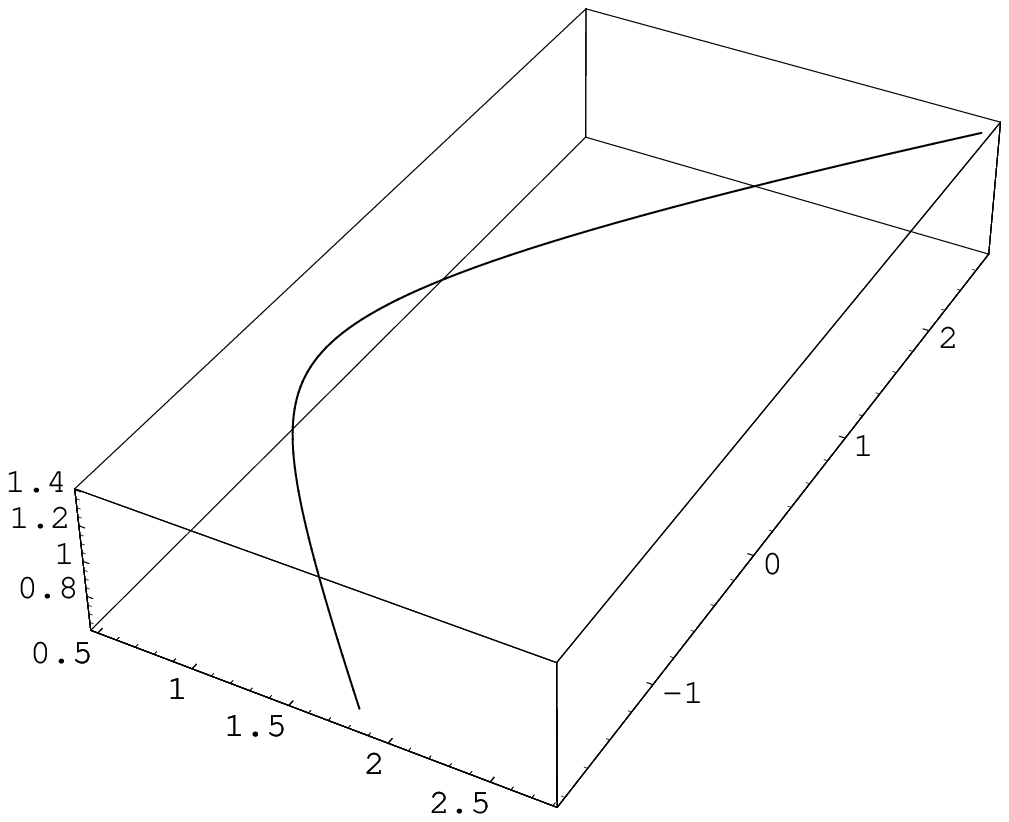}
\hspace*{0.5cm}
\includegraphics[width=4.0cm]{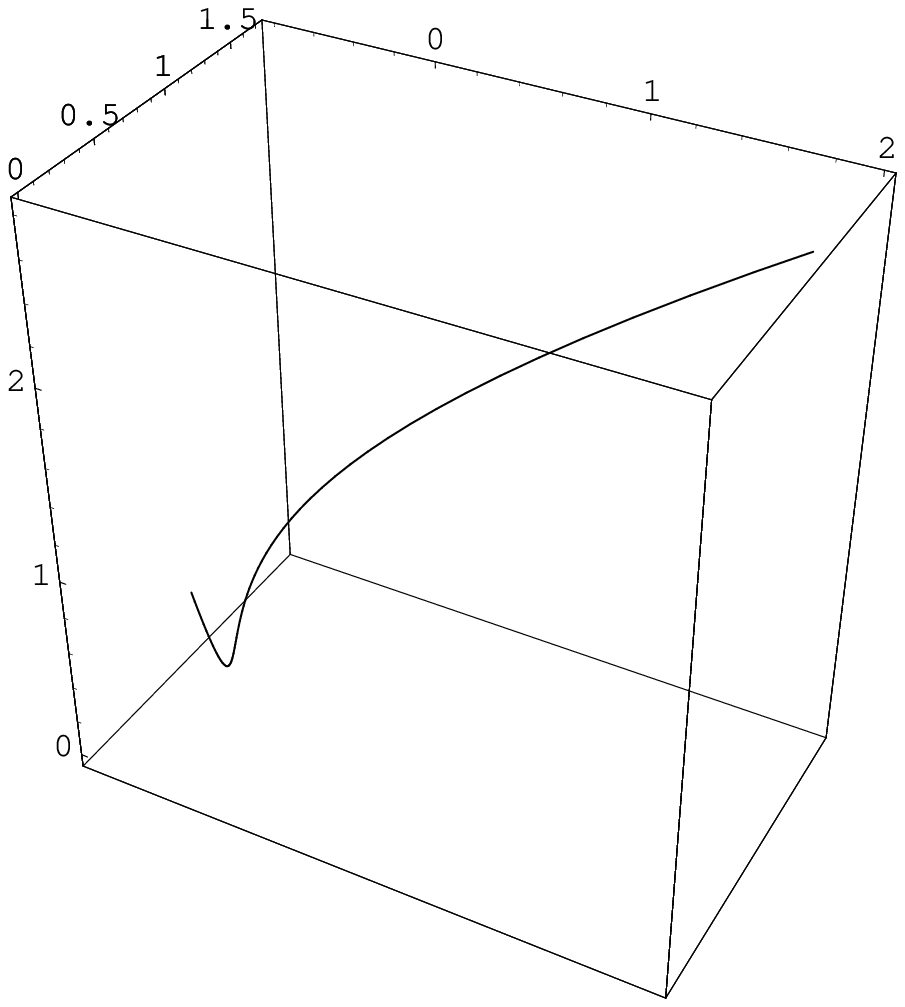}
\hspace*{0.5cm}
\includegraphics[width=4.0cm]{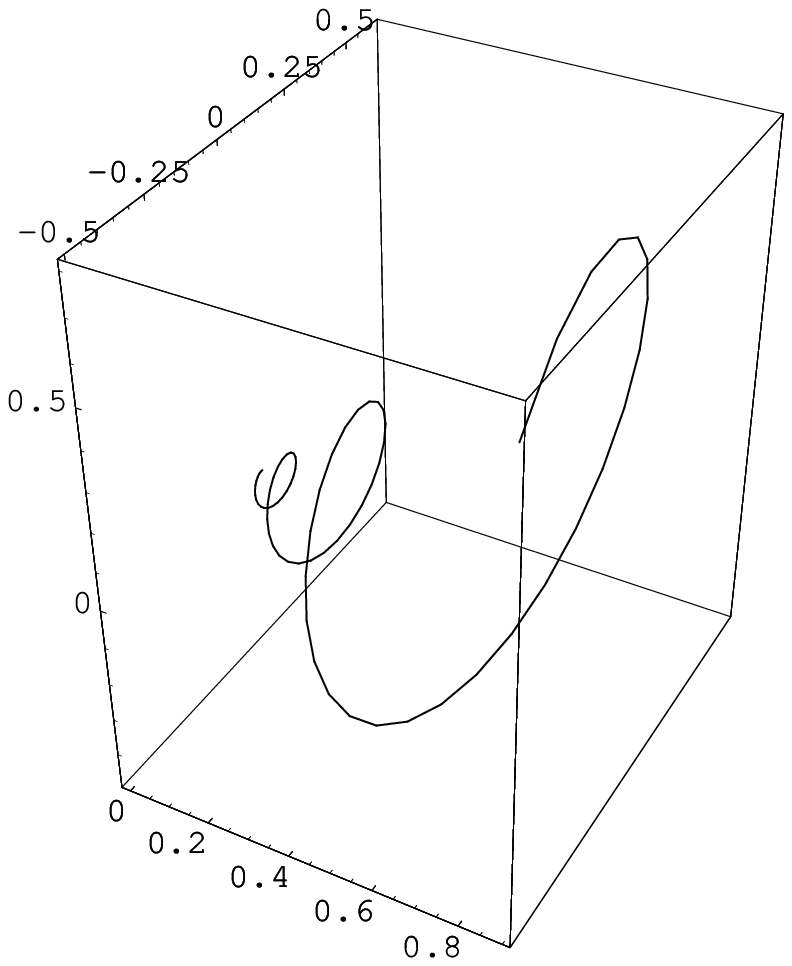}
\end{center}\caption{Some spacelike general helix with $\kappa=\frac{h}{s}$ and $\tau=\frac{r}{s}$.}
\label{fig-x3}
\end{figure}

\end{document}